\documentclass{article}

\newcommand{\sfZF}{{\sf ZF}}

\newcommand{\fsubset}{\subset_{fin}}
\newcommand{\dep}{\mbox{{\rm dp}}}

\usepackage{proof}
\usepackage{latexsym}
\usepackage{amsfonts}
\usepackage{amssymb}
\usepackage{cite}

\newcommand{\blem}{\begin{lemma}}
\newcommand{\elem}{\end{lemma}}
\newcommand{\bth}{\begin{theorem}}
\newcommand{\ethm}{\end{theorem}}
\newcommand{\benu}{\begin{enumerate}}
\newcommand{\eenu}{\end{enumerate}}
\newcommand{\bdes}{\begin{description}}
\newcommand{\edes}{\end{description}}
\newcommand{\bdf}{\begin{definition}}
\newcommand{\edf}{\end{definition}}
\newcommand{\bcor}{\begin{cor}}
\newcommand{\ecor}{\end{cor}}
\newcommand{\bprp}{\begin{proposition}}
\newcommand{\eprp}{\end{proposition}}
\newcommand{\bmlem}{\begin{mlemma}}
\newcommand{\emlem}{\end{mlemma}}
\newcommand{\bclm}{\begin{claim}}
\newcommand{\eclm}{\end{claim}}
\newcommand{\bprf}{{\bf Proof}.\hspace{2mm}}

\newcommand{\eprf}{\hspace*{\fill} $\Box$}

\newcommand{\beqn}{\begin{equation}}
\newcommand{\eeqn}{\end{equation}}
\newcommand{\beqnarr}{\begin{eqnarray}}
\newcommand{\eeqnarr}{\end{eqnarray}}
\newcommand{\beqnarrs}{\begin{eqnarray*}}
\newcommand{\eeqnarrs}{\end{eqnarray*}}

\newcommand{\spand}{\,\&\,}

\newtheorem{theorem}{Theorem}[section]
\newtheorem{definition}[theorem]{Definition}
\newtheorem{proposition}[theorem]{Proposition}
\newtheorem{lemma}[theorem]{Lemma}
\newtheorem{cor}[theorem]{Corollary}
\newtheorem{mlemma}[theorem]{Main Lemma}
\newtheorem{claim}[theorem]{Claim}

\newcommand{\alp}{\alpha}

\newcommand{\veps}{\varepsilon}

\newcommand{\Del}{\Delta}
\newcommand{\ome}{\omega}
\newcommand{\Ome}{\Omega}
\newcommand{\bet}{\beta}

\newcommand{\Gam}{\Gamma}
\newcommand{\kap}{\kappa}

\newcommand{\Sig}{\Sigma}
\newcommand{\tht}{\theta}
\newcommand{\Tht}{\Theta}

\newcommand{\vphi}{\varphi}

\newcommand{\fal}{\forall}
\newcommand{\exi}{\exists}

\newcommand{\Rarw }{\Rightarrow}

\newcommand{\lrarw}{\leftrightarrow}
\newcommand{\Lrarw}{\Leftrightarrow}

\newcommand{\calh}{{\cal H}}

\newcommand{\calk}{{\cal K}}

\newcommand{\sfk}{{\sf k}}

\newcommand{\la}{\langle}
\newcommand{\ra}{\rangle}
\newcommand{\lc}{\lceil}
\newcommand{\rc}{\rceil}

\title{
Conservations of first-order reflections
}
\author{Toshiyasu Arai
\\
Graduate School of Science,
Chiba University
\\
1-33, Yayoi-cho, Inage-ku,
Chiba, 263-8522, JAPAN
\\
tosarai@faculty.chiba-u.jp
}
\date{}
\begin{document}
\maketitle

\begin{abstract}
The set theory {\sf KP}$\Pi_{N+1}$ for $\Pi_{N+1}$-reflecting universes is shown to be
$\Pi_{N+1}$-conservative over iterations of 
$\Pi_{N}$-recursively Mahlo operations for each $N\geq 2$.
\end{abstract}


\section{Introduction}\label{sect:introduction}

It is well known that the set of weakly Mahlo cardinals below a weakly compact cardinal
is stationary.
Furthermore any weakly compact cardinal $\kap$ is in the diagonal intersection $\kap\in M^{\triangle}=\bigcap\{M(M^{\alp}):\alp<\kap\}$
for the $\alp$-th iterate $M^{\alp}$ of the Mahlo operation $M$, 
where $\kap\in M(X)$ iff $X\cap\kap$ is stationary in $\kap$.

The same holds for the recursive analogues of the indescribable cardinals, \textit{reflecting ordinals}
introduced by Richter-Aczel\cite{Richter-Aczel74}.
First let us recall the ordinals briefly.
For a full account of the admissible set theory, see \cite{barwise}.

$\Del_{0}$ denotes the set of bounded formulas in the language $\{\in\}$ of set theories.
Then the classes $\Sig_{i+1},\Pi_{i+1}$ are defined recursively as usual.
Each class $\Sig_{i+1},\Pi_{i+1}$ is defined to be closed under bounded quantifications $\exi x\in a,\fal x\in a$.

The axioms of the Kripke-Platek set theory with the axiom of infinity, denoted $\mbox{{\sf KP}}\ome$,
are Extensionality, Foundation schema, Pair, Union, $\Del_{0}$-Separation, $\Del_{0}$-Collection, and the axiom of infinity.
Note that except Foundation schema,
each axiom in {\sf KP}$\ome$ is a $\Pi_{2}$-formula.

For set-theoretic formulas $\varphi$, let
$
P\models\varphi:\Leftrightarrow (P,\in)\models\varphi
$.

In what follows, let $V$ denote a transitive and wellfounded model of $\mbox{{\sf KP}}\ome$, which is a universe in discourse.
$P, Q,\ldots$ denote non-empty transitive sets in $V\cup\{V\}$.

A $\Pi_{i}$-recursively Mahlo operation for $2\leq i<\omega$, is defined through a universal $\Pi_{i}$-formula 
$\Pi_{i}(a)$: 
\begin{eqnarray*}
P\in RM_{i}({\mathcal X}) & :\Leftrightarrow & \forall b\in P[P\models\Pi_{i}(b) \to \exists Q\in {\mathcal X}\cap P(b\in Q\models\Pi_{i}(b))]
\\
&&\mbox{(read:} P \mbox{ is } \Pi_{i}\mbox{-reflecting on } {\mathcal X}\mbox{.)}
\end{eqnarray*}
For the universe $V$, $V\in RM_{i}(\mathcal{X})$ denotes 
$\forall b[\Pi_{i}(b) \to \exists Q\in {\mathcal X}(b\in Q\models\Pi_{i}(b))]$.
Suppose that there exists a first-order sentence $\varphi$ such that $P\in {\mathcal X} \Leftrightarrow P \models \varphi$ for any transitive $P\in V\cup\{V\}$.
Then  $RM_{i}({\mathcal X})$ is $\Pi_{i+1}$, i.e., there exists a $\Pi_{i+1}$-sentence $rm_{i}({\mathcal X})$
such that $P\in RM_{i}({\mathcal X})$ iff $P\models rm_{i}({\mathcal X})$ for any transitive 
set $P$.

The iteration of $RM_{i}$ along a definable relation $\prec$ is defined as follows.
\[
P\in RM_{i}(a;\prec) :\Lrarw a\in P\in\bigcap\{ RM_{i}(RM_{i}(b;\prec)): b\in P\models b\prec a\}.
\]
Again $P\in RM_{i}(a;\prec)$ is a $\Pi_{i+1}$-relation.

Let $Ord$ denote the class of ordinals in $V$.
Let us write $RM_{i}^{\alp}$ for $RM_{i}(\alp;<)$ and ordinals $\alp\in Ord$.
A transitive set $P$ is said to be $\Pi_{i}${\it -reflecting\/} 
if $P\in RM_{i}=RM_{i}^{1}$.

$P\in RM_{i+1}$ is much stronger than $P\in RM_{i}$: Assume $P\in RM_{i+1}$ and 
$P\models\Pi_{i}(b)$ for $b\in P$.
Then $P\in RM_{i}$ and $P\models rm_{i}\land\Pi_{i}(b)$ for the $\Pi_{i+1}$-sentence $rm_{i}$
such that $P\in RM_{i}$ iff $P\models rm_{i}$.
Hence there exists a $Q\in P$ such that $Q\models rm_{i}\land\Pi_{i}(b)$, i.e.,
$Q\in RM_{i} \,\&\, Q\models\Pi_{i}(b)$.
This means $P\in RM_{i}^{2}=RM_{i}(RM_{i})$.
Moreover $P$ is in the diagonal intersection of $RM_{i}$,
$P\in RM_{i}^{\triangle}$, i.e., $P\in\bigcap\{RM_{i}^{\beta}:\beta\in P\cap Ord\}$, and so on.

In particular a set theory KP$\Pi_{i+1}$  for universes in $RM_{i+1}$ proves the consistency of a set theory 
for universes in $RM_{i}^{\triangle}$.

In this paper we address a problem:
How far can we iterate lower recursively Mahlo operations in higher reflecting universes?
In \cite{LMPS} we gave a sketchy proof of the following Theorem \ref{th:resolution0}, which is implicit in
 ordinal analyses in \cite{WienpiN, ptpiN}.

\bth\label{th:resolution0}
For each $N\geq 2$ there exists a $\Sig_{1}$-relation $\lhd_{N}$ on $\ome$
such that
the set theory {\sf KP}$\ell$ for limits of admissibles proves the transfinite induction schema for $\lhd_{N}$
up to each $a\in\omega$, and
{\sf KP}$\Pi_{N+1}$ is $\Pi^{1}_{1}$(on $\ome$)-conservative over
the theory 
$${\sf KP}\ell+\{V\in RM_{N}(a;\lhd_{N}): a\in\ome\}.$$
\end{theorem}

Theorem \ref{th:resolution0} suffices to approximate {\sf KP}$\Pi_{N+1}$ proof-theoretically
in terms of iterations of 
$\Pi_{N}$-recursively Mahlo operations.
However $V\in RM_{N}(a;\prec)$ is a $\Pi_{N+1}$-formula for $\Sig_{N+1}$-relation $\prec$,
and the class $\Pi^{1}_{1}$ on $\ome$ is smaller than $\Pi_{N+1}$.

In this paper the set theory {\sf KP}$\Pi_{N+1}$ for $\Pi_{N+1}$-reflecting universes is shown to be
$\Pi_{N+1}$-conservative over iterations of 
$\Pi_{N}$-recursively Mahlo operations $RM_{N}$ for each $N\geq 2$.
This result will be extended in \cite{liftupK, liftupKn} to the indescribable cardinals over 
$\sfZF+(V=L)$.

\section{Conservation}\label{sect:Consv}

Let $Ord\subset V$ denote the class of ordinals, $Ord^{\veps}\subset V$ and $<^{\veps}$ be $\Del$-predicates such that
for any transitive and wellfounded model $V$ of $\mbox{{\sf KP}}\ome$,
$<^{\veps}$ is a well ordering of type $\veps_{\Ome+1}$ on $Ord^{\veps}$
for the order type $\Ome$ of the class $Ord$ in $V$.
Specifically let us encode `ordinals' $\alp<\veps_{\Ome+1}$ by codes $\lceil\alp\rceil\in Ord^{\veps}$ as follows.
$\lceil\alp\rceil=\la 0,\alp\ra$ for $\alp\in Ord$, 
$\lceil \Ome\rceil=\la 1,0\ra$,
$\lc\ome^{\alp}\rc=\la 2,\lc\alp\rc\ra$ for $\alp>\Ome$,
and
$\lc \alp\rc=\la 3,\lc\alp_{1}\rc,\ldots,\lc\alp_{n}\rc\ra$
if $\alp=\alp_{1}+\cdots+\alp_{n}>\Ome$ with $\alp_{1}\geq\cdots\geq\alp_{n}$, $n>1$ and 
$\exi\bet_{i}(\alp_{i}=\ome^{\bet_{i}})$ for each $\alp_{i}$.
Then
$\lceil\ome_{n}(\Ome+1)\rceil\in Ord^{\veps}$ denotes the code of the `ordinal' 
$\ome_{n}(\Ome+1)$.

$<^{\veps}$ is assumed to be a canonical ordering such that 
$\mbox{{\sf KP}}\ome$ proves the fact that $<^{\veps}$ is a linear ordering, and for any formula $\vphi$
and each $n<\ome$,
\beqn\label{eq:trindveps}
\mbox{{\sf KP}}\ome\vdash\fal x(\fal y<^{\veps}x\,\vphi(y)\to\vphi(x)) \to \fal x<^{\veps}\lceil\ome_{n}(\Ome+1)\rceil\vphi(x)
\eeqn
For a definition of $\Del$-predicates $Ord^{\veps}$ and $<^{\veps}$, 
and a proof of (\ref{eq:trindveps}), cf. \cite{liftupZF}.

\begin{proposition}\label{prp:transitive}
${\sf KP}\ome$ proves that if
$P\in RM_{N}(\bet;<^{\veps})$, then $ \fal\alp<^{\veps}\bet(\alp\in P \to P\in RM_{N}(\alp;<^{\veps}))$.
\end{proposition}
\bprf
This is seen from the fact that $<^{\veps}$ is transitive in {\sf KP}$\ome$.
\eprf

\bth\label{th:resolutionpi3}

For each $N\geq 2$,
{\sf KP}$\Pi_{N+1}$ is $\Pi_{N+1}$-conservative over
the theory 
\[
\mbox{{\sf KP}}\ome+
\{
V\in RM_{N}(\lceil\ome_{n}(\Ome+1)\rceil;<^{\veps}): 
n\in\ome\}
.\]
\end{theorem}

From (\ref{eq:trindveps}) we see that {\sf KP}$\Pi_{N+1}$ proves $V\in RM_{N}(\lceil\ome_{n}(\Ome+1)\rceil;<^{\veps})$
for each $n\in\ome$.
The converse is proved in section \ref{sect:proof}.

\begin{proposition}
For any class $\Gam$ of $\Pi_{N+1}$-sentences, there exists a $\Sig_{N+1}$-sentence $A$ such that
${\sf KP}\Pi_{N+1}\vdash A$, and ${\sf KP}\ome+\Gam\not\vdash A$ unless ${\sf KP}\ome+\Gam$ is inconsistent.
\end{proposition}
\bprf
This follows from the essential unboundedness theorem due to G. Kreisel and A. L\'evy \cite{KL}.
In this proof let $\vdash A:\Lrarw {\sf KP}\ome\vdash A$ and ${\rm Pr}$ denote a standard provability predicate for 
{\sf KP}$\ome$.
Also ${\rm Tr}_{\Pi_{N+1}}$ denotes a partial truth definition of $\Pi_{N+1}$-sentences.

Then let $A$ be a $\Sig_{N+1}$-sentence saying that `I am not provable from any true $\Pi_{N+1}$-sentence',
$\vdash A\lrarw \fal x\in\ome[{\rm Tr}_{\Pi_{N+1}}(x) \to \lnot{\rm Pr}(x\dot{\to}\lc A\rc)]$,
where $\dot{\to}$ denotes a recursive function such that
$\lc A\rc\dot{\to}\lc B\rc=\lc A\to B\rc$ for codes
$\lc A\rc$ of formulas $A$.

Suppose ${\sf KP}\ome+\Gam\vdash A$. Pick a $C\in\Gam$ so that $\vdash C\to A$.
Then ${\sf KP}\ome+\Gam\vdash{\rm Tr}_{\Pi_{N+1}}(\lc C\rc)\land {\rm Pr}(\lc C\to A\rc)$.
Hence ${\sf KP}\ome+\Gam\vdash \lnot A$.

In what follows argue in ${\sf KP}\Pi_{N+1}$.
Suppose $A$ is false, and let $C$ be any true $\Pi_{N+1}$-sentence.
Since the universe $V$ is $\Pi_{N+1}$-reflecting,
there exists a transitive model $P\in V$ of ${\sf KP}\ome+\{C,\lnot A\}$, which shows that
${\sf KP}\ome+\{C,\lnot A\}$ is consistent.
In other words, $\lnot{\rm Pr}(\lc C\to A\rc)$.
Therefore ${\sf KP}\Pi_{N+1}\vdash \lnot A\to A$.
\eprf
\\

Thus Theorem \ref{th:resolutionpi3} is optimal with respect to the class $\Pi_{N+1}$ of formulas
provided that {\sf KP}$\Pi_{N+1}$ is consistent.

\bcor
For each $N\geq 3$,
{\sf KP}$\Pi_{N+1}+(\mbox{{\rm Power}})+(\Sig_{N-2}\mbox{{\rm -Separation}})+(\Pi_{N-2}\mbox{{\rm -Collection}})$ is $\Pi_{N+1}$-conservative over
the theory 
$
\mbox{{\sf KP}}\ome+
\{
V\in RM_{N}(\lceil\ome_{n}(\Ome+1)\rceil;<^{\veps}): 
n\in\ome\}+(\mbox{{\rm Power}})+(\Sig_{N-2}\mbox{{\rm -Separation}})+(\Pi_{N-2}\mbox{{\rm -Collection}})
$.
\ecor
\bprf
This follows from Theorem \ref{th:resolutionpi3} and 
the facts that the axiom Power is a $\Pi_{3}$-sentence $\fal a\exi b\fal x\subset a(x\in b)$, and 
$\Sig_{i}$-Separation or $\Pi_{i}$-Collection are $\Pi_{i+2}$-formulas.
\eprf
\\

Let us announce an extension of Theorem \ref{th:resolutionpi3} in \cite{liftupK, liftupKn} to the indescribable cardinals over $\sfZF+(V=L)$.

Let $<^{\veps}$ be an $\varepsilon$-ordering as above.
Let $M_{N}$ denote the $\Pi^{1}_{N}$-Mahlo operation defined for sets $S$ of ordinals and uncountable regular cardinals $\kap$:
$\kap\in M_{N}(S)$ iff $S\cap\kap$ is $\Pi^{1}_{N}$-indescribable in $\kap$.
The $\Pi^{1}_{N+1}$-indescribability is proof-theoretically reducible
to iterations of an operation along initial segments of $<^{\veps}$ over $\sfZF+(V=L)$.
The operation is a mixture of the $\Pi^{1}_{N}$-Mahlo operation $M_{N}$ and Mostowski collapsings.

For $\alp<^{\veps}\veps_{\calk+1}$ and finite sets $\Tht\fsubset(\calk+1)$,
$\Pi_{n+1}$-classes $Mh_{n}^{\alp}[\Tht]$ are defined so that the following holds.

In Theorem\ref{th:mainthKN} $\calk$ is intended to denote the least $\Pi^{1}_{N+1}$-indescribable cardinal, 
and $\Ome$ the least weakly inaccessible cardinal above $\calk$.

\bth\label{th:mainthKN}{\rm (The case $N=0$ in \cite{liftupK}, and the general case in \cite{liftupKn}.)}

\benu
\item\label{th:mainthK1}
{\rm For each} $n<\ome$,
\[
\sfZF+(V=L)+
(\calk \mbox{ {\rm is }} \Pi^{1}_{N+1}\mbox{{\rm -indescribable}})\vdash \calk\in Mh_{n}^{\ome_{n}(\Ome+1)}[\emptyset]
.\]

\item\label{th:mainthK2}
{\rm For any} $\Sig^{1}_{N+2}${\rm -sentences} $\vphi$, 
{\rm if}
\[
\sfZF+(V=L)+
(\calk \mbox{ {\rm is }} \Pi^{1}_{N+1}\mbox{{\rm -indescribable}})\vdash\vphi^{L_{\calk}}
,
\]
{\rm then we can find an} $n<\ome$ {\rm such that}
\[
\sfZF+(V=L)+
(\calk\in Mh_{n}^{\ome_{n}(\Ome+1)}[\emptyset])
\vdash
\vphi^{L_{\calk}}
.\]

\eenu

\end{theorem}

The classes $Mh_{n}^{\alp}[\Tht]$ are defined from iterated Skolem hulls $\calh_{\alp,n}(X)$, through which 
we described the limit of $\sfZF+(V=L)$-provable countable ordinals in \cite{liftupZF} as follows.

\bth\label{th:mainthZ}{\rm (\cite{liftupZF}.)}
\beqnarrs
|\sfZF+(V=L)|_{\ome_{1}} & := &
 \inf\{\alp\leq\ome_{1}: \fal \vphi [\sfZF+(V=L) \vdash \exi x\in L_{\ome_{1}}\,\vphi
\Rarw \exi x\in L_{\alp}\, \vphi]\}
\\
& = &
\Psi_{\ome_{1}}\veps_{\Ome+1}:=\sup\{\Psi_{\ome_{1},n}\ome_{n}(\Ome+1): n<\ome\}.
\eeqnarrs

\end{theorem}
In Theorem \ref{th:mainthZ}, $\Ome$ is intended to denote the least weakly inaccessible cardinal.

\section{Proof of Theorem \ref{th:resolutionpi3}}\label{sect:proof}

In this section we prove Theorem \ref{th:resolutionpi3}.
Our proof is extracted from M. Rathjen's ordinal analyses of $\Pi_{3}$-reflection in \cite{Rathjen94}.

Let $N\geq 2$ denote a fixed integer.
The axioms of the set theory {\sf KP}$\Pi_{N+1}$ for $\Pi_{N+1}$-reflecting universes
are those of {\sf KP}$\ome$
and the axiom for $\Pi_{N+1}$-reflection: for $\Pi_{N+1}$-formulas $\vphi$,
$\vphi(a) \to \exi c[ad^{c}\land a\in c \land \vphi^{c}(a)]$,
where $ad$ denotes a $\Pi_{3}$-sentence such that $P\models ad$ iff $P$ is a transitive model of {\sf KP}$\ome$,
and
$\vphi^{c}$ denotes the result of restricting any unbounded quantifiers $\exi x,\fal x$ in $\vphi$ to $\exi x\in c,\fal x\in c$, resp.

$\mbox{{\sf KP}}i$ denotes the set theory for recursively inaccessible sets, which is obtained
from {\sf KP}$\ome$ by adding the axiom
$\fal x\exi y[x\in y\land ad^{y}]$.

Throughout this section we work in an intuitionistic fixed point theory 
$\mbox{FiX}^{i}(\mbox{{\sf KP}}i)$ over $\mbox{{\sf KP}}i$.
The intuitionistic theory $\mbox{FiX}^{i}(\mbox{{\sf KP}}i)$ is introduced in \cite{liftupZF},
and shown to be a conservative extension of $\mbox{{\sf KP}}i$.
Let us reproduce definitions and results on $\mbox{FiX}^{i}(\mbox{{\sf KP}}i)$ here.

Fix an $X$-strictly positive formula $\mathcal{Q}(X,x)$ in the language $\{\in,=,X\}$ with an extra unary predicate symbol $X$.
In $\mathcal{Q}(X,x)$ the predicate symbol $X$ occurs only strictly positive.
This means that the predicate symbol $X$ does not occur in the antecedent $\vphi$ of implications $\vphi\to\psi$ 
nor in the scope of negations $\lnot$ in $\mathcal{Q}(X,x)$.
The language of $\mbox{FiX}^{i}(\mbox{{\sf KP}}i)$ is $\{\in,=,Q\}$ with a fresh unary predicate symbol $Q$.
The axioms in $\mbox{FiX}^{i}(\mbox{{\sf KP}}i)$ consist of the following:
\benu
\item
All provable sentences in $\mbox{{\sf KP}}i$ (in the language $\{\in,=\}$).

\item
Induction schema for any formula $\vphi$ in $\{\in,=,Q\}$:
\beqn\label{eq:Qind}
\fal x(\fal y\in x\,\vphi(y)\to\vphi(x))\to\fal x\,\vphi(x)
\eeqn

\item
Fixed point axiom:
\[
\fal x[Q(x)\lrarw \mathcal{Q}(Q,x)]
.\]
\eenu

The underlying logic in $\mbox{FiX}^{i}(\mbox{{\sf KP}}i)$ is defined to be the intuitionistic (first-order predicate) logic (with equality).

(\ref{eq:Qind}) yields the following Lemma \ref{lem:vepsfix}.

\blem\label{lem:vepsfix}
Let $<^{\veps}$ denote a $\Del_{1}$-predicate mentioned in the beginning of section \ref{sect:Consv}.
For each $n<\ome$ and each formula $\vphi$ in $\{\in,=,Q\}$,
\[
\mbox{{\rm FiX}}^{i}(\mbox{{\sf KP}}i)\vdash\fal x(\fal y<^{\veps}x\,\vphi(y)\to\vphi(x)) \to 
\fal x<^{\veps}\ome_{n}(I+1)\vphi(x)
.\]
\elem

The following Theorem \ref{th:consvintfix} is seen as in \cite{intfix, liftupZF}.

\bth\label{th:consvintfix}
$\mbox{{\rm FiX}}^{i}(\mbox{{\sf KP}}i)$ is a conservative extension of $\mbox{{\sf KP}}i$.
\end{theorem}

In what follows we work in $\mbox{FiX}^{i}(\mbox{{\sf KP}}i)$.

Let $V$ denote a transitive and wellfounded model of $\mbox{{\sf KP}}\ome$.
Consider the language $\mathcal{L}_{V}=\{\in\}\cup\{c_{a}:a\in V\}$ where $c_{a}$ denotes the name of the set $a\in V$.
We identify the set $a$ with its name $c_{a}$.

Our proof proceeds as follows.
Assume that {\sf KP}$\Pi_{N+1}\vdash A$ for a $\Pi_{N+1}$-sentence $A$.
{\sf KP}$\Pi_{N+1}$ is embedded to an infinitary system $R_{N}$ formulated in one-sided sequent calculus,
and cut inferences are eliminated, which results in an infinitary derivation of height $\alp<\veps_{\Ome+1}$
with an inference rule $(Ref_{N+1})$ for $\Pi_{N+1}$-reflection.
Then $A$ is seen to be true in $P\in RM_{N}(\alp;<^{\veps})$.

In one-sided sequent calculi, formulas are generated from atomic formulas and their negations
$a\in b,a\not\in b$
by propositional connectives $\lor,\land$ and quantifiers $\exi,\fal$.
It is convenient here to have bounded quantifications $\exi x\in a,\fal x\in a$ besides unbounded ones $\exi x,\fal x$.
The negation $\lnot A$ of formulas $A$ is defined recursively by de Morgan's law and elimination of double negations.
Also $(A\to B):\equiv(\lnot A\lor B)$.

$\Gam,\Del,\ldots$ denote finite sets of sentences, called \textit{sequents} in the language $\mathcal{L}_{V}$.
$\Gam,\Del$ denotes the union $\Gam\cup\Del$, and $\Gam,A$ the union $\Gam\cup\{A\}$.
A finite set $\Gam$ of sentences is intended to denote the disjunction $\bigvee\Gam:=\bigvee\{A:A\in\Gam\}$.
$\Gam$ is \textit{true} in $P\in V\cup\{V\}$ 
iff $\bigvee\Gam$ is true in $P$ iff $\bigvee\Gam^{P}$ is true.

Classes $\Del_{0},\Sig_{i+1},\Pi_{i+1}$ of sentences in $\mathcal{L}_{V}$ are defined as usual.

We assign disjunctions or conjunctions to sentences as follows.
When a disjunction $\bigvee(A_{\iota})_{\iota\in J}$ [a conjunction $\bigwedge(A_{\iota})_{\iota\in J}$]
is assigned to $A$,
we denote $A\simeq\bigvee(A_{\iota})_{\iota\in J}$ [$A\simeq\bigwedge(A_{\iota})_{\iota\in J}$], resp.

\bdf\label{df:assigndc}
\benu
\item
{\rm For a} $\Del_{0}${\rm -sentence} $M$
\[
M:\simeq
\left\{
\begin{array}{ll}
\bigvee(A_{\iota})_{\iota\in J} & \mbox{{\rm if }} M \mbox{ {\rm is false in }} V
\\
\bigwedge(A_{\iota})_{\iota\in J} &  \mbox{{\rm if }} M \mbox{ {\rm is true in }} V
\end{array}
\right.
\mbox{{\rm with }} J:=\emptyset
.\]

{\rm In what follows consider the unbounded sentences.}

\item
$(A_{0}\lor A_{1}):\simeq\bigvee(A_{\iota})_{\iota\in J}$
{\rm and}
$(A_{0}\land A_{1}):\simeq\bigwedge(A_{\iota})_{\iota\in J}$
{\rm with} $J:=2$.

\item
$
\exi x\in a\, A(x):\simeq\bigvee(A(b))_{b\in J}$
{\rm and}
$
\fal x\in a\, A(x):\simeq\bigwedge(A(b))_{b\in J}
$ 
{\rm with}
$
J:=a$.

\item
$
\exi x\, A(x):\simeq\bigvee(A(b))_{b\in J}
$ 
{\rm and}
$
\fal x\, A(x):\simeq\bigwedge(A(b))_{b\in J}
$ 
{\rm with}
$
J:=V
$.

\eenu
\edf

\bdf
{\rm The} \textit{depth} $\dep(A)<\ome$ {\rm of} $\mathcal{L}_{V}${\rm -sentences} $A$ {\rm is defined recursively as follows.}
\benu
\item
$\dep(A)=0$ {\rm if} $A\in\Del_{0}$.

{\rm In what follows consider unbounded sentences} $A$.

\item
$\dep(A)=\max\{\dep(A_{i}):i<2\}+1$ {\rm if} $A\equiv (A_{0}\circ A_{1})$ 
{\rm for} $\circ\in\{\lor,\land\}$.
\item
$\dep(A)=\dep(B(\emptyset))+1$ {\rm if} $A\in\{ (Q x \, B(x)), (Q x\in a\, B(x)) :a\in V\}$
for $Q\in\{\exi,\fal\}$.
\eenu
\edf

\bdf
\benu
\item
{\rm For} $\mathcal{L}_{V}${\rm -sentences} $A$,
$
\sfk(A):=\{a\in V: c_{a}\mbox{ {\rm occurs in} } A\}
$.

\item
{\rm For sets} $\Gam$ {\rm of sentences,} $\sfk(\Gam):=\bigcup\{\sfk(A):A\in\Gam\}$.


\item
{\rm For} $\iota\in V$ {\rm and a transitive model}
 $P\in V\cup\{V\}$ {\rm of} {\sf KP}$\ome$,
$P(\iota)\in V\cup\{V\}$ 
{\rm denotes the smallest transitive model of}
 {\sf KP}$\ome$ {\rm such that} $P\cup\{\iota\}\subset P(\iota)$, {\rm cf. \cite{barwise}.}

{\rm For finite lists} $\vec{a}=(a_{1},\ldots,a_{n})$, $P(\vec{a}):=(\cdots P(a_{1})\cdots)(a_{n})$.

\eenu
\edf




Inspired by operator controlled derivations due to W. Buchholz \cite{Buchholz},
let us define a relation $P\vdash^{\alp}_{m}\Gam$ for transitive models $P\in V\cup\{V\}$ of {\sf KP}$\ome$.
The relation $P\vdash^{\alp}_{m}\Gam$ is defined as a fixed point of a strictly positive formula $H$
\[
H(P,\alp,m,\Gam)\Lrarw P\vdash^{\alp}_{m}\Gam
\]
in $\mbox{FiX}^{i}(\mbox{{\sf KP}}i)$.

Note that $P$ contains the code $\la 1,0\ra=\lc\Ome\rc$, and is closed under ordinal addition $(\alp,\bet)\mapsto\alp+\bet$, exponentiation $\alp\mapsto\ome^{\alp}$ for $\alp,\bet\in Ord^{\veps}$
and $a\mapsto rank(a)$ for $rank(a)=\sup\{rank(b)+1:b\in a\}$.


\bdf\label{df:controlderreg}
{\rm Let} $P\in V\cup\{V\}$ {\rm be a transitive model of} {\sf KP}$\ome$,
$\alp<\veps_{\Ome+1}$ {\rm and} $m<\ome$.

$P\vdash^{\alp}_{m}\Gam$ {\rm holds if}
\beqn\label{eq:controlord}
\sfk(\Gam)\cup\{\alp\}\subset P
\eeqn
{\rm and one of the following
cases holds:}

\bdes
\item[$(\bigvee)$]
{\rm there is an} $A\in\Gam$ {\rm such that}
$A\simeq\bigvee(A_{\iota})_{\iota\in J}$, {\rm and for an} $\iota\in J$  {\rm and an}  $\alp(\iota)<\alp$,
$P\vdash^{\alp(\iota)}_{m}\Gam,A_{\iota}$.
\[
\infer[(\bigvee)]{P\vdash^{\alp}_{m}\Gam}{P\vdash^{\alp(\iota)}_{m}\Gam,A_{\iota}}
\]

\item[$(\bigwedge)$]
{\rm there is an} $A\in\Gam$ {\rm such that}
$A\simeq\bigwedge(A_{\iota})_{\iota\in J}$, {\rm and for any} $\iota\in J$, 
{\rm there is an} $\alp(\iota)$ {\rm such that} $\alp(\iota)<\alp$ {\rm and}
$P(\iota)\vdash^{\alp(\iota)}_{m}\Gam,A_{\iota}$.
\[
\infer[(\bigwedge)]{P\vdash^{\alp}_{m}\Gam}{\{P(\iota)\vdash^{\alp(\iota)}_{m}\Gam,A_{\iota}:\iota\in J\}}
\]

\item[$(cut)$]
{\rm there are} $C$ {\rm and} $\alp_{0},\alp_{1}$ {\rm such that}
$\dep(C)<m$, $\alp_{0},\alp_{1}<\alp$, {\rm and} $P\vdash^{\alp_{0}}_{m}\Gam,\lnot C$ {\rm and}
 $P\vdash^{\alp_{1}}_{m}C,\Gam$.

\[
\infer[(cut)]{P\vdash^{\alp}_{m}\Gam}
{P\vdash^{\alp_{0}}_{m}\Gam,\lnot C & P\vdash^{\alp_{1}}_{m}C,\Gam}
\]

\item[$(Ref_{N+1})$]
{\rm there are} $A(c)\in\Pi_{N+1}$ {\rm and} 
$\alp_{0},\alp_{1}<\alp$ {\rm such that}
$P\vdash^{\alp_{0}}_{m}\Gam,A(c)$
{\rm and} 
$P\vdash^{\alp_{1}}_{m}\fal z[ad^{z} \to c\in z \to \lnot A^{z}],\Gam$.

\[
\infer[(Ref_{N+1})]{P\vdash^{\alp}_{m}\Gam}
{
P\vdash^{\alp_{0}}_{m}\Gam,A(c)
&
P\vdash^{\alp_{1}}_{m}\fal z[ad^{z} \to c\in z \to \lnot A^{z}],\Gam
}
\]
\edes
\edf


In what follows, let us fix an integer $n_{0}$ and restrict (codes of) ordinals to 
$\alp<^{\veps}\lc\ome_{n_{0}}(\Ome+1)\rc$.
$n_{0}$ is chosen from the given finite proof of a $\Pi_{N+1}$-sentence $A$ in {\sf KP}$\Pi_{N+1}$, cf.
Corollary \ref{th:embedreg} (Embedding).
Since $n_{0}$ is a constant, we see from Lemma \ref{lem:vepsfix} that
$\mbox{FiX}^{i}(\mbox{{\sf KP}}i)$ proves
 transfinite induction schema up to $\lc\ome_{n_{0}}(\Ome+1)\rc$
 for any formula
in which the derivability relation $P\vdash^{\alp}_{m}\Gam$ may occur.

\begin{proposition}\label{lem:vdash}
Let $P^{\prime}\supset P$ be a transitive model of ${\sf KP}\ome$,  $m\leq m^{\prime}<\ome$
and 
$\sfk(\Del)\cup\{\alp^{\prime}\}\subset P^{\prime}$.
If $P\vdash^{\alp}_{m}\Gam$,  then 
$P^{\prime}\vdash^{\alp^{\prime}}_{m^{\prime}}\Gam,\Del$.
\end{proposition}

In embedding ${\sf KP}\Pi_{N+1}$ in the infinitary calculus,
it is convenient to formulate ${\sf KP}\Pi_{N+1}$ in (finitary) one-sided sequent calculus
of the language $\{\in,0\}$ with the individual constant $0$ for the empty set.
Axioms are logical ones $\Gam,\lnot A,A$ for any formulas $A$,
and axioms in the theory ${\sf KP}\Pi_{N+1}$.
Inference rules are $(\lor)$, $(\land)$ for propositional connectives, $(b\exi)$, $(b\fal)$ for
bounded quantifications, $(\exi)$, $(\fal)$ for unbounded quantifications,
and $(cut)$. 
For details, see the proof of the following Lemma \ref{th:embedregax}.

Though the following Lemmas \ref{th:embedregax}, \ref{lem:reduction} and \ref{lem:predcereg} are seen as in \cite{Buchholz},
we give proofs of them for readers' convenience.

Let 
$(m,\vec{a}):=\Ome\cdot m+3 rank(a_{1})\#\cdots\# 3 rank(a_{n})$
for $\vec{a}=(a_{1},\ldots,a_{n})$ and the natural (commutative) sum $\alp\#\bet$ of ordinals $\alp,\bet$.

\blem\label{th:embedregax}
Suppose ${\sf KP}\Pi_{N+1}\vdash\Gam(\vec{x})$, 
where free variables occurring in the sequent are among the list $\vec{x}$.
Then there is an $m<\ome$ 
such that for any $\vec{a}\subset V$ and any transitive model $P\in V\cup\{V\}$ of {\sf KP}$\ome$,
$P(\vec{a})\vdash_{m}^{(m,\vec{a})}\Gam(\vec{a})$.
\elem
\bprf
First consider the logical axiom $\Gam(\vec{x}),\lnot A(\vec{x}),A(\vec{x})$.
We see that for any $\vec{a}$
\beqn\label{eq:embedtaut}
P(\vec{a})\vdash^{2d}_{0}\Gam(\vec{a}),\lnot A(\vec{a}),A(\vec{a})
\eeqn
by induction on $d=\dep(A)$.

Then by Proposition \ref{lem:vdash} we have $P(\vec{a})\vdash^{(2d,\vec{a})}_{2d}\Gam(\vec{a}),\lnot A(\vec{a}),A(\vec{a})$.

If $d=0$, then $A\in\Del_{0}$ and one of $\lnot A(\vec{a})$ and $A(\vec{a})$ is true.
Hence by $(\bigwedge)$ we have $P(\vec{a})\vdash^{0}_{0}\Gam(\vec{a}),\lnot A(\vec{a}),A(\vec{a})$.

Next consider the case when $A\equiv (\exi y\, B(\vec{x},y))\not\in\Del_{0}$ with $\dep(B(\vec{x},y))=d-1$.
By IH(=Induction Hypothesis) we have for any $\vec{a}\subset V$ and any $b\in V$,
$P(\vec{a}*(b))\vdash^{2d-2}_{0}\Gam(\vec{a}),\lnot B(\vec{a},b),B(\vec{a},b)$, where
$(a_{1},\ldots,a_{n})*(b)=(a_{1},\ldots,a_{n},b)$.
$(\bigvee)$ yields $P(\vec{a}*(b))\vdash^{2d-1}_{0}\Gam(\vec{a}),\lnot B(\vec{a},b),\exi y\, B(\vec{a},y)$.
Hence $(\bigwedge)$ with $P(\vec{a}*(b))=P(\vec{a})(b)$ yields $P(\vec{a})\vdash^{2d}_{0}\Gam(\vec{a}),\lnot \exi y\, B(\vec{a},y),\exi y\, B(\vec{a},y)$.

The cases $A\equiv (\exi y\in a\, B(\vec{x},y))\not\in\Del_{0}$ and $A\equiv(B_{0}\lor B_{1})\not\in\Del_{0}$ are similar.
Thus  (\ref{eq:embedtaut}) was shown.

Second consider the inference rule $(\exi)$ with $\exi y\, A(\vec{x},y)\in\Gam(\vec{x})$
\[
\infer[(\exi)]{\Gam(\vec{x})}{\Gam(\vec{x}),A(\vec{x},t)}
\]
When $t$ is a variable $y$, we can assume that $y$ is an $x_{i}$ in the list $\vec{x}$,
for otherwise substitute $0$ for $y$.
By IH there is an $m$ such that
$P(\vec{a})\vdash_{m}^{(m,\vec{a})}\Gam(\vec{a}),A(\vec{a},t^{\prime})$
where $t^{\prime}\equiv a_{i}$ if $t\equiv x_{i}$, and $t^{\prime}\equiv 0$ otherwise.
Thus $P(\vec{a})\vdash_{m+1}^{(m+1,\vec{a})}\Gam(\vec{a})$.

Third consider the inference rule $(\fal)$ with $\fal y\, A(\vec{x},y)\in\Gam(\vec{x})$
\[
\infer[(\fal)]{\Gam(\vec{x})}{\Gam(\vec{x}),A(\vec{x},y)}
\]
where the variable $y$ does not occur in $\Gam(\vec{x})$.
IH yields for an $m$,
$P(\vec{a}*(b))\vdash_{m}^{(m,\vec{a}*(b))}\Gam(\vec{a}),A(\vec{a},b)$.
$(\bigwedge)$ with $(m+1,\vec{a})>(m,\vec{a}*(b))$ yields $P(\vec{a})\vdash_{m+1}^{(m+1,\vec{a})}\Gam(\vec{a})$.

The following cases are similarly seen.
\[
\infer[(b\exi)]{\Gam,\exi y\in s\, B(\vec{x},y)}{\Gam,t\in s & \Gam, B(\vec{x},t)},
\:
\infer[(b\fal)]{\Gam,\fal y\in s\, B(\vec{x},y)}{\Gam,y\not\in s, B(\vec{x},y)},
\:
\infer[(\lor)]{\Gam,A_{0}\lor A_{1}}{\Gam,A_{0},A_{1}},
\:
\infer[(\land)]{\Gam, A_{0}\land A_{1}}{\Gam,A_{0} & \Gam,A_{1}}
\]

In a cut inference
\[
\infer[(cut)]{\Gam(\vec{x})}{\Gam(\vec{x}),\lnot A(\vec{x}) & A(\vec{x}),\Gam(\vec{x})}
\]
if the cut formula $A(\vec{x})$ has free variables $\vec{y}$ other than $\vec{x}$, then substitute $0$ for $\vec{y}$.

In what follows let us suppress parameters.

Fourth consider the axioms other than Foundation.
For example, consider the $\Del_{0}$-Collection $\fal x\in a\exi y\, A(x,y)\to \exi z\fal x\in a\exi y\in z\, A(x,y)$
for $A\in\Del_{0}$ and $a\in V$.
Since $P(a)$ is a transitive model of {\sf KP}$\ome$ and $a\in P(a)$,
pick a $b\in P(a)$ such that $\fal x\in a\exi y\, A(x,y)\to \fal x\in a\exi y\in b\, A(x,y)$ holds in $P(a)$.
Then $\lnot\fal x\in a\exi y\, A(x,y)\lor\fal x\in a\exi y\in b\, A(x,y)$ is a true $\Del_{0}$-sentence.
Hence
$P(a)\vdash_{0}^{0}\lnot\fal x\in a\exi y\, A(x,y),\fal x\in a\exi y\in b\, A(x,y)$.
Three $(\bigvee)$'s yield $P(a)\vdash_{0}^{3}\fal x\in a\exi y\, A(x,y)\to \exi z\fal x\in a\exi y\in z\, A(x,y)$.

Next consider the axiom $A(c)\to\exi z[ad^{z}\land c\in z\land A^{z}]$ for $A\in\Pi_{N+1}$.
We have by (\ref{eq:embedtaut}) for $d=\dep(A)$
\[
\infer[(Ref_{N+1})]
{P(c)\vdash_{0}^{2d+1}\lnot A(c),\exi z[ad^{z}\land c\in z\land A^{z}]}
{P(c)\vdash^{2d}_{0}\lnot A(c),A(c) & P(c)\vdash_{0}^{2}\fal z[ad^{z}\to c\in z\to \lnot A^{z}],\exi z[ad^{z}\land c\in z\land A^{z}]}
\]

In this way we see that there are cut-free infinitary derivations of finite heights deducing axioms in {\sf KP}$\Pi_{N+1}$
other than Foundation.

Finally consider the Foundation.
Let $d=\dep(A)$ and $B\equiv(\lnot\fal x(\fal y\in x\, A(y)\to A(x)))$.
We show by induction on $rank(a)$ that 
\beqn\label{eq:embedfund}
P(a)\vdash_{0}^{2d+3 rank(a)}B, \fal x\in a\, A(x)
\eeqn
By IH we have for any $b\in a$,
$P(b)\vdash_{0}^{2d+3 rank(b)}B, \fal x\in b\, A(x)$.
Thus we have by (\ref{eq:embedtaut})
\[
\infer[(\bigvee)]{P(b)\vdash_{0}^{2d+3 rank(b)+2}B, A(b)}
{
\infer[(\bigwedge)]{P(b)\vdash_{0}^{2d+3 rank(b)+1}B, \fal x\in b\, A(x)\land \lnot A(b),A(b)}
{
 P(b)\vdash_{0}^{2d+3 rank(b)}B, \fal x\in b\, A(x)
 &
 P(b)\vdash_{0}^{2d}\lnot A(b),A(b)
 }
}
\]
Therefore (\ref{eq:embedfund}) was shown.
\[
\infer[(\bigwedge)]
{P(a)\vdash_{0}^{2d+3 rank(a)}B,\fal x\in a\, A(x)}
{\{P(a,b)\vdash_{0}^{2d+3 rank(b)+2}B, A(b): b\in a\}}
\]
 
\eprf

\bcor\label{th:embedreg}{\rm (Embedding)}
If ${\sf KP}\Pi_{N+1}\vdash A$ for a sentence $A$, then there is an $m<\ome$ such that for any transitive model $P\in V\cup\{V\}$ of {\sf KP}$\ome$,
 $P\vdash_{m}^{\Ome\cdot m}A$.
\ecor

\blem\label{lem:reduction}{\rm (Reduction)}
Let $C\simeq \bigvee(C_{\iota})_{\iota\in J}$. Then
\[
(P\vdash^{\alp}_{m}\Del,\lnot C) \spand (P\vdash^{\bet}_{m}C,\Gam) \spand (\dep(C)\leq m)
\Rarw
P\vdash^{\alp+\bet}_{m}\Del,\Gam
\]
\elem
\bprf
This is seen by induction on $\bet$.

Consider first the case when $C$ is a $\Del_{0}$-sentence.
Then $C$ is false and $J=\emptyset$.
From $P\vdash^{\bet}_{m}C,\Gam$ we see that $P\vdash^{\bet}_{m}\Gam$.
$\bet\leq\alp+\bet$ yields $P\vdash^{\alp+\bet}_{m}\Del,\Gam$.

Next assume that the last inference rule in $P\vdash^{\bet}_{m}C,\Gam$ is a $(\bigvee)$ with the main formula 
$C\not\in\Del_{0}$:
\[
\infer[(\bigvee)]{P\vdash^{\bet}_{m}C,\Gam}
{P\vdash^{\bet(\iota)}_{m}C,C_{\iota},\Gam}
\]
where $\iota\in J$ and $\bet(\iota)<\bet$.
We can assume that $\iota$ occurs in $C_{\iota}$.
Otherwise set $\iota=0$.
Thus
$\iota\in P$ by (\ref{eq:controlord}).
On the other hand we have $P(\iota)\vdash^{\alp}_{m}\Del,\lnot C_{\iota}$ by inversion,
and hence $P\vdash^{\alp}_{m}\Del,\lnot C_{\iota}$ by $\iota\in P$.

IH yields
$P\vdash^{\alp+\bet(\iota)}_{m}C_{\iota},\Del,\Gam$.
A cut inference with $P\vdash^{\alp}_{m}\Del,\lnot C_{\iota}$ and $\dep(C_{\iota})<\dep(C)\leq m$ 
yields $P\vdash^{\alp+\bet}_{m}\Del,\Gam$.
\eprf

\blem\label{lem:predcereg}{\rm (Predicative Cut-elimination)}
$P\vdash^{\alp}_{m+1}\Gam
 \Rarw 
P\vdash^{\ome^{\alp}}_{m}\Gam$.
\elem
\bprf
This is seen by induction on $\alp$ using Reduction \ref{lem:reduction} and the fact: $\bet<\alp\Rarw \ome^{\bet}+\ome^{\bet}\leq\ome^{\alp}$.
\eprf
\\

For $\alp<^{\veps}\lc\ome_{n}(\Ome+1)\rc$, set $RM_{N}^{\alp}:=RM_{N}(\alp;<^{\veps})$.

\begin{proposition}\label{prp:Mahloness}
Let $\Gam\subset\Pi_{N+1}$ and $P\in RM_{N}^{\alp}$ be a transitive model of {\sf KP}$i$.
Assume
\[
\exi\xi,x\in P(0<^{\veps}\xi<^{\veps}\alp\land
\fal Q\in RM_{N}^{\xi}\cap P(x\in Q\models{\sf KP}i \to \Gam \mbox{ {\rm is true in }} Q))
.\]
Then
$\Gam$ is true in $P$.
\end{proposition}
\bprf
By $P\in RM_{N}^{\alp}$ we have
$P\in RM_{N}(RM_{N}^{\xi})$ for any $\xi\in P$ such that $\xi<^{\veps}\alp$.

Suppose the $\Sig_{N+1}$-sentence $\vphi:=\bigwedge\lnot\Gam:=\bigwedge\{\lnot\tht:\tht\in\Gam\}$ is true in $P$.
Then for any $\xi\in P$ with $\xi<^{\veps}\alp$ and $x\in P$ there exists a transitive model $Q\in RM_{N}^{\xi}\cap P$
of {\sf KP}$i$ minus $\Del_{0}$-Collection
such that $x\in Q$ and $\vphi$ is true in $Q$.
Let $0<^{\veps}\xi<^{\veps}\alp$.
Then any $Q\in RM_{N}^{\xi}$ is $\Pi_{N}$-reflecting for $N\geq 2$, and hence
 is a model of $\Del_{0}$-Collection.
\eprf

\blem\label{lem:CollapsingthmKR}{\rm (Elimination of  $(Ref_{N+1})$)}
Let $\Gam\subset\Pi_{N+1}$.
Suppose $P_{0}\vdash^{\alp}_{0}\Gam$, $P_{0}\in P$ and $P\in RM_{N}^{\alp}$ for a transitive model $P$ of {\sf KP}$i$.
Then $\Gam$ is true in $P$.
\elem
\bprf This is seen by induction on $\alp$. 
Let $P_{0}\vdash^{\alp}_{0}\Gam$, $P_{0}\in P$ and $P\in RM_{N}^{\alp}$ be a transitive model $P$ of {\sf KP}$i$.
Note that any sentence occurring in the witnessed derivation of $P_{0}\vdash^{\alp}_{0}\Gam$
is $\Pi_{N+1}$.
\\

\noindent
{\bf Case 1}.
First consider the case when the last inference is a $(Ref_{N+1})$:
By (\ref{eq:controlord}) we have $\{\alp_{\ell},\alp_{r}\}\subset P_{0}\subset P$, $\max\{\alp_{\ell},\alp_{r}\}<^{\veps}\alp$, $A\in\Pi_{N+1}$.
\[
\infer[(Ref_{N+1})]{P_{0}\vdash^{\alp}_{0}\Gam}
{
P_{0}\vdash^{\alp_{\ell}}_{0}\Gam,A(c)
&
P_{0}\vdash^{\alp_{r}}_{0}\fal z[ad^{z} \to c\in z \to \lnot A^{z}(c)],\Gam
}
\]
If $\alp_{\ell}=0$, then there is a $B\in\Gam\cup A(c)$ such that $B\simeq\bigwedge(B_{\iota})_{\iota\in\emptyset}$, i.e.,
$B$ is either a true $\Del_{0}$-sentence or a sentence $\fal x\in\emptyset\, C(x)$.
In each case we can assume $B\in\Gam$, and $B^{P}\in\Gam^{P}$ is true.

In what follows assume $0<^{\veps}\alp_{\ell}$.
We can assume that $c$ occurs in $A(c)$, and hence 
$c\in P_{0}$.

By Proposition \ref{prp:transitive} we have $P\in RM_{N}^{\alp_{r}}$. 
From IH we see that 
\beqn\label{eq:case1r}
\mbox{either } \fal z\in P[ad^{z} \to c\in z \to \lnot A^{z}(c)] 
\mbox{ or } \bigvee\Gam^{P} \mbox{ is true.}
\eeqn

On the other hand by IH we have for any $Q\in RM_{N}^{\alp_{\ell}}\cap P$ with $c\in P_{0}\in Q\models{\sf KP}i$ 
that
either $\bigvee\Gam^{Q}$ is true or $A(c)^{Q}$ is true.
By (\ref{eq:case1r})
for any $Q\in RM_{N}^{\alp_{\ell}}\cap P$ with $P_{0}\in Q\models{\sf KP}i$,
$\bigvee\Gam^{Q}\lor\bigvee\Gam^{P}$ is true.
From Proposition \ref{prp:Mahloness} with $0<^{\veps}\alp_{\ell}$ we see that $\bigvee\Gam^{P}$ is true.
\\

\noindent
{\bf Case 2}.
Second consider the case when the last inference is a $(\bigwedge)$:
we have $A\simeq\bigwedge(A_{\iota})_{\iota\in J}$, $A\in\Gam$, and $\alp(\iota)<\alp$ for any $\iota\in J$
\[
\infer[(\bigwedge)]{P_{0}\vdash^{\alp}_{0}\Gam}{\{P_{0}(\iota)\vdash^{\alp(\iota)}_{0}\Gam,A_{\iota}:\iota\in J\}}
\]
For any $\iota\in P$ we have $P_{0}(\iota)\in P$ since $P$ is assumed to be a limit of admissibles.

IH yields for any $\iota\in P$ that
either $\bigvee\Gam^{P}$ is true or $A_{\iota}^{P}$ is true.
If $J=V$, then we are done.
If $J=a\in V$, then $a\in P_{0}\subset P$ by (\ref{eq:controlord}) , and hence $a\subset P$.
\\

\noindent
{\bf Case 3}.
Third consider the case when the last inference is a $(\bigvee)$:
we have $A\simeq\bigvee(A_{\iota})_{\iota\in J}$, $A\in\Gam$, and $\alp(\iota)<\alp$ for an $\iota\in J$
\[
\infer[(\bigvee)]{P_{0}\vdash^{\alp}_{0}\Gam}{P_{0}\vdash^{\alp(\iota)}_{0}\Gam,A_{\iota}}
\]
IH yields that
either $\bigvee\Gam^{P}$ is true or $A_{\iota}^{P}$ is true.
Consider the case when $J=V$.
We can assume that $\iota$ occurs in $A_{\iota}$.
Then $\iota\in P_{0}\subset P$.
Hence $\bigvee\Gam^{P}$ is true.
\eprf
\\

Let us prove Theorem \ref{th:resolutionpi3}.
Let $N\geq 2$, and $A$ be a $\Pi_{N+1}$-sentence provable in {\sf KP}$\Pi_{N+1}$.
By Embedding \ref{th:embedreg} and Predicative Cut-elimination \ref{lem:predcereg} we have for an $m<\ome$,
$L_{\ome_{1}^{CK}}\vdash_{0}^{\ome_{m}(\Ome\cdot m)}A$ with $\emptyset=\sfk(A)$, and 
$L_{\ome_{1}^{CK}}\vdash_{0}^{\ome_{n}(\Ome+1)}A$ for an $n>m$ with $n<n_{0}$.
If $V\in RM_{N}^{\ome_{n}(\Ome+1)}$, then $L_{\ome_{1}^{CK}}\in V\models{\sf KP}i$, and $A$ is true (in $V$)
by Elimination of  $(Ref_{N+1})$ \ref{lem:CollapsingthmKR}.

By formalizing the above proof in $\mbox{FiX}^{i}(\mbox{{\sf KP}}i)$ with 
Lemma \ref{lem:vepsfix}
yields 
\[
\mbox{FiX}^{i}(\mbox{{\sf KP}}i)\vdash V\in RM_{N}(\lceil\ome_{n}(\Ome+1)\rceil;<^{\veps})\to A
\]
and then by Theorem \ref{th:consvintfix} and 
$\mbox{{\sf KP}}\ome\vdash V\in RM_{N}(\lceil\ome_{n}(\Ome+1)\rceil;<^{\veps})\to
\fal x\exi y[x\in y\land ad^{y}]$
\[
\mbox{{\sf KP}}\ome\vdash V\in RM_{N}(\lceil\ome_{n}(\Ome+1)\rceil;<^{\veps})\to A
.\]

In the formalization note that we have in $\mbox{FiX}^{i}(\mbox{{\sf KP}}i)$, 
a partial truth definition of $\Pi_{N+1}$-sentences.

\bibliographystyle{amsplain}

\end{document}